# Robust Risk-Constrained Unit Commitment with Large-scale Wind Generation: An Adjustable Uncertainty Set Approach

Cheng Wang, Feng Liu, *Member, IEEE*, Jianhui Wang, *Senior Member*, *IEEE*, Feng Qiu, *Member, IEEE,* Wei Wei, *Member*, *IEEE*, Shengwei Mei, *Fellow*, *IEEE*, Shunbo Lei, *Student Member*, *IEEE*

*Abstract—* This paper addresses two vital issues which are barely discussed in the literature on robust unit commitment (RUC): 1) how much the potential operational loss could be if the realization of uncertainty is beyond the prescribed uncertainty set; 2) how large the prescribed uncertainty set should be when it is used for RUC decision making. In this regard, a robust risk-constrained unit commitment (RRUC) formulation is proposed to cope with large-scale volatile and uncertain wind generation. Differing from existing RUC formulations, the wind generation uncertainty set in RRUC is adjustable via choosing diverse levels of operational risk. By optimizing the uncertainty set, RRUC can allocate operational flexibility of power systems over spatial and temporal domains optimally, reducing operational cost in a risk-constrained manner. Moreover, since impact of wind generation realization out of the prescribed uncertainty set on operational risk is taken into account, RRUC outperforms RUC in the case of rare events. The traditional column and constraint generation (C&CG) and two algorithms based on C&CG are adopted to solve the RRUC. As the proposed algorithms are quite general, they can also apply to other RUC models to improve their computational efficiency. Simulations on a modified IEEE 118-bus system demonstrate the effectiveness and efficiency of the proposed methodology.

## NOMENCLATURE

*Indices*

| | |
|---|---|
| $g$ | Index for generators. |
| $m$ | Index for wind farms. |
| $l$ | Index for transmission lines. |
| $j$ | Index for loads. |
| $n$ | Index for nodes. |
| $t$ | Index for time periods. |

*Parameters*

| | |
|---|---|
| $T$ | Number of time periods. |
| $N$ | Number of nodes. |
| $M$ | Number of wind farms. |
| $G$ | Number of generators. |
| $L$ | Number of transmission lines. |
| $S_g$ | Start-up cost for generator $g$. |
| $c_g$ | Constant term of generation cost function. |
| $P_g^{min}/P_g^{max}$ | Minimal/ maximal output of generator $g$. |
| $R_+^g/R_-^g$ | Ramp-up/ ramp-down limit for generator $g$. |
| $T_g^{on}/T_g^{off}$ | Minimum on/off hour of generator $g$. |
| $F_l$ | Transmission capacity of line $l$. |
| $W$ | Wind generation uncertainty set. |
| $\widehat{w}_{mt}$ | Forecasted output of wind farm $m$ in period $t$. |
| $w_m^{max}$ | Installed capacity of wind farm $m$. |
| $\Gamma^S/\Gamma^T$ | Uncertainty budget over spatial/ temporal scale. |
| $D_{jt}$ | Load demand of load node $j$ in period $t$. |
| $B$ | Node admittance matrix of the grid. |
| $Line_l$ | Indices of initial node and terminal node of line $l$. |
| $o_1/o_2$ | Number of initial/ terminal node of line $l$. |
| $\Phi(n)$ | The set of nodes connecting to node $n$. |
| $\alpha_{mt}$ | Confidence level of wind generation output interval of wind farm $m$ in period $t$. |
| $\beta_t/\beta_s$ | Confidence level of $\Gamma^T/\Gamma^S$. |
| $e_t$ | Price of wind generation curtailment in period $t$. |
| $f_t$ | Price of load shedding in period $t$. |
| $\pi_{gt}/\pi_{mt}/\pi_{jt}$ | Generation shift distribution factor of generator $g$/ wind farm $m$/ load $j$ in period $t$. |
| $Risk_{dh}$ | Day-ahead operational risk level. |

*Decision Variables*

| | |
|---|---|
| $u_{gt}$ | Binary variable indicating whether generator $g$ is on or off in period $t$. |
| $z_{gt}$ | Binary variable indicating whether generator $g$ is started up in period $t$. |
| $p_{gt}$ | Real-time output of generator $g$ in period $t$. |
| $\hat{p}_{gt}$ | Day-ahead output of generator $g$ in period $t$. |
| $v_{mt}^u/v_{mt}^l$ | Binary variable indicating normalized positive /negative output deviation of wind farm $m$ in period $t$. |
| $\Delta w_{mt}$ | Wind generation curtailment in wind farm $m$ in |

This work was supported in part by China State Grid Corp Science and Technology Project (SGSXDKY-DWKJ2015-001), Foundation for Innovative Research Groups of the National Natural Science Foundation of China (51321005), the Special Fund of National Basic Research Program of China (2012CB215103), and the fund of China Scholarship Council (CSC).

C. Wang, F. Liu, W. Wei, and S. Mei are with the State Key Laboratory of Power Systems, Department of Electrical Engineering and Applied Electronic Technology, Tsinghua University, 100084 Beijing, China. (E-mail: c-w12@mails.tsinghua.edu.cn; lfeng@mail.tsinghua.edu.cn; wei-wei04@mails.tsinghua.edu.cn; meishengwei@tsinghua.edu.cn).

J. Wang and F. Qiu are with the Argonne National Laboratory, Argonne, IL 60439, USA (e-mail: jianhui.wang@anl.gov; fqiu@anl.gov).

S. Lei is with Department of Electrical and Electronic Engineering, The University of Hong Kong, Pokfulam, Hong Kong (email: leishunbo@eee.hku.hk).

| | |
|---|---|
| | period $t$. |
| $\Delta D_{jt}$ | Load shedding at load node $j$ in period $t$. |
| $w^u_{mt}$/ $w^l_{mt}$ | Upper/ lower bound of wind generation output of wind farm $m$ in period $t$. |
| $Q^p_{mt}$ /$Q^n_{mt}$ | Operational risk due to underestimation/ overestimation of the output of wind farm $m$ in period $t$. |
| $\theta_{nt}$ | Phase angle of node $n$ in period $t$. |

## I. INTRODUCTION

THE increasing penetration of wind generation has brought great challenges to power system operation and scheduling. In day-ahead scheduling, the main challenge is how to make reliable and economic dispatch decisions to effectively hedge against considerable volatile and uncertainty wind generation. Unit commitment (UC) decision making is among the most important issues as it essentially determines the operational flexibility of a power system in the following day. Many inspiring works on this topic have been done, which can be roughly divided into two categories: scenario based stochastic unit commitment (SUC) and uncertainty set based robust unit commitment (RUC).

In [1], [2], [3], several stochastic unit commitment (SUC) models have been proposed. SUC can generate an economic and reliable strategy over a set of selected scenarios. However, it may miss out certain scenarios with a small chance of occurrence albeit a severe adverse consequence. Thus the resulting strategy may be vulnerable to those rare unfavorable scenarios. In this context, robust unit commitment (RUC) is more reliable as it can guarantee the operational feasibility for every possible scenario in the prescribed uncertainty set [4], [5]. This approach, however, may increase conservativeness of UC strategy in return. To circumvent this problem, many models and methods are developed, such as the minimax regret unit commitment [6], the unified stochastic and robust unit commitment [7], the hybrid stochastic/interval approach [8], and the multi-band uncertainty set approach [9]. Similar strategy conservatism issues also exist in robust economic dispatch (RED) problems. In [10], dynamic uncertainty sets considering the temporal and spatial relationship correlation of uncertainty are adopted in RED. In [11], the authors shrink the uncertainty bands to balance dispatch cost and dispatch infeasibility penalty based on a given prediction interval.

Most existing research on robust dispatch, including RUC and RED, assume that the wind generation uncertainty set is given, within which the operational feasibility of power system can be completely guaranteed. However, one crucial issue is barely discussed in the literature: how large the potential loss could be if the realizations of uncertain wind generation are beyond the scope of the prescribed uncertainty set. According to [12], in United States, the real-time wind generation can deviate even more than 6 and 10 times the standard deviation from day-ahead and hour-ahead forecast values, respectively. In practice, reserve adequacy assessment (RAA) or reliability assessment commitment (RAC) would be carried out between clearing of day-ahead (DA) market and real-time (RT) market, which would improve the capability of guaranteeing the balance of demand and supply [13], [14]. Such rare events may still cause operational infeasibility as well as operational loss even though its probability may be small, as the penetration of wind generation keeps increaseing. This necessitates a risk measure to quantify the potential operational loss due to the rear events beyond the prescribed uncertainty set and facilitate the UC decision making.

The above problem consequently raises another important issue: how large the prescribed uncertainty set should be when it is utilized for robust UC decision making. To address the aforementioned concern, a novel concept called do-not-exceed (DNE) limit is proposed in [15], in which the admissible wind generation interval for each wind farm can be obtained under a fixed economic dispatch (ED) strategy. However, this approach does not take fully advantage of wind generation forecast information and the selection of coefficients in its objective function is a hot potato. Another extensively adopted treatment is to use probability density function (PDF) of wind generation as well as a unified confidence level to determine the parameters of the uncertainty set [5], [16], [17]. This treatment, however, may encounter two major obstacles in practice: 1) the determination of confidence level is subjective to a large extent; 2) there is lack of systematic methods to determine an appropriate confidence level for each of dispatch periods. In [18], appropriate uncertainty set can be obtained from simulation experiments, which avoid the subjectivity of choosing confidence level, yet may be time-consuming as the number of time period and uncertainty source increase.

In this paper, the expectation of operational loss is used as a risk measure to depict the impact of operational infeasibility on operation as well as a metric for system operator to adjust the range of the uncertainty set. This risk measure considers not only the consequences of potential operational infeasibility, but also the probability of those consequences. In the literature, much effort has been devoted to incorporate risk management into power system dispatch. In [19], the inspiring concept of risk-limiting dispatch (RLD) is proposed, where the risk is calculated under a given acceptable loss of load probability (LOLP) and load shedding (LS) cost coefficient. [20] presents a risk-based UC model for day-ahead market clearing. [21] suggests a risk constrained robust unit commitment model in which uncertainty sets of multiple sources are considered. [22] demonstrates a chance-constrained unit commitment model with *n-k* security criterion in which the conditional value-at-risk (CVaR) is minimized. In [23], the CVaR-based transaction cost is minimized in a robust optimal power flow formulation.

This paper derives a novel risk measure, which is referred to as *operational risk*, includes expected operational loss for wind generation curtailment (WGC) as well as LS. Based on this, a robust risk-constrained UC model is formulated by taking the operational risk as a constraint rather than as an objective function in conventional RUC models. Column and constraints generation (C&CG) based algorithms are also proposed to solve the model. Compared with existing works, major contributions of this paper are summarized in twofold.

1) The mathematical formulation of RRUC is proposed, which is formulated as a two-stage robust optimization model. In the first stage, the sum of UC cost and ED cost under predicted value of wind generation are minimized with an explicit constraint of operational risk. In the second stage, the feasibility of the first-stage decision variables against wind generation uncertainty is guaranteed. The salient features of the RRUC are listed as follows.

i. The boundaries of wind generation uncertainty set are first-stage adjustable decision variables in RRUC rather than parameters in RUC. By optimizing the boundaries of uncertainty set, RRUC allows an optimal allocation of the operational flexibility of the power system over spatial and temporal domains.

ii. Based on our previous work [24], the operational risk is utilized as an additional constraint added into RRUC model, enabling it to strictly guarantee the operational feasibility in normal scenarios within the uncertainty set while limiting the operational risk under rare events beyond the uncertainty set.

iii. Due to the discontinuity of UC strategy, the admissible boundaries of wind generation under RUC decisions are always larger than the prescribed uncertainty set. To obtain exact admissible boundaries of uncertain wind generation, risk-based admissibility assessment [24] or DNE limit [15] assessment has to be conducted after the UC decision is given. In RRUC, nevertheless, the optimal solution concurrently gives the exact admissible boundaries of wind generation without additional computation.

2) Mathematically, the RRUC formulation leads to a two-stage robust optimization model. Algorithms for RUC such as such as Benders decomposition (BD) and C&CG algorithms can be directly applied. However, due to the modeling difference, efficiency of these algorithms may be undermined. In this paper, two C&CG-based algorithms are developed to reduce the number of iterations as well as computational scale. They could also be applied to other two-stage robust optimization models.

The remaining part of the paper is organized as follows. Section II presents the mathematical formulation. Section III derives the solution methodology. Section IV gives the case studies to demonstrate the proposed model and algorithms. Finally, section V concludes the paper with discussion.

## II. MATHEMATICAL FORMULATION

### A. Operational Risk under Wind Generation Uncertainty

The concept of wind generation admissibility region (WGAR) is proposed in [24] (See Fig. 1). [24] also reveals that the boundaries of WGAR can be obtained by solving a wind generation admissibility assessment problem. If the realization of wind generation is within WGAR, there will be no any operational loss in the following day (in other words, WGAR is riskless). The rest part of the region of wind generation is defined as the inadmissibility region (WGIR). If any realization of uncertain wind generation intersects with WGIR, certain operational loss may occur in the scheduling day. The operational risk in WGIR can be defined as

$$Risk = \sum_{t=1}^{T}\sum_{m=1}^{M} \left( e_t \int_{w_{mt}^u - \hat{w}_{mt}}^{w_m^{max} - \hat{w}_{mt}} \left( \delta_{mt} - w_{mt}^u + \hat{w}_{mt} \right) + \cdots \right. \\ \left. + f_t \int_{-\hat{w}_{mt}}^{w_{mt}^l - \hat{w}_{mt}} \left( w_{mt}^l - \delta_{mt} - \hat{w}_{mt} \right) \right) y_{mt}(\delta_{mt}) \, d\delta_{mt} \quad (1a)$$

where, $w_{mt}^u$ and $w_{mt}^l$ are the upper and the lower boundaries of admissible wind generation, respectively. $e_t$ and $f_t$ are cost coefficients of WGC and LS, respectively. $\delta_{mt}$ represents the wind generation forecast error and $y_{mt}(\cdot)$ is its PDF. In (1a), the first and second integral terms represent the operational risk caused by underestimated and overestimated wind generation, respectively. Formula (1a) can be approximated by a linear expression with auxiliary variables and constraints using piecewise linearization (PWL) method as follows.

$$Risk = \min_{Q_{mt}^p, Q_{mt}^n} \sum_{t=1}^{T}\sum_{m=1}^{M}(Q_{mt}^p + Q_{mt}^n) \quad (1b)$$

$$Q_{mt}^p \geq a_{mtsz}^p w_{mt}^u + b_{mtsz}^p \quad \forall m, \forall t, s=0,1\ldots,S, z=0,1,\ldots,Z-1. \quad (1c)$$

$$Q_{mt}^n \geq a_{mtsz}^n w_{mt}^l + b_{mtsz}^n \quad \forall m, \forall t, s=0,1\ldots,S, z=0,1,\ldots,Z-1. \quad (1d)$$

where, (1b) is the linear approximation of (1a); (1c) and (1d) are the auxiliary constraints induced by the PWL treatment. $a_{mtsz}^p, a_{mtsz}^n, b_{mtsz}^p, b_{mtsz}^n$ are constant coefficients of the piecewise linear approximation; $s$ and $z$ are ordinal numbers generated during the PLA treatment; $S$ and $Z$ are the maximum values of $s$ and $z$, respectively. We refer the readers to [24] for more details on model (1).

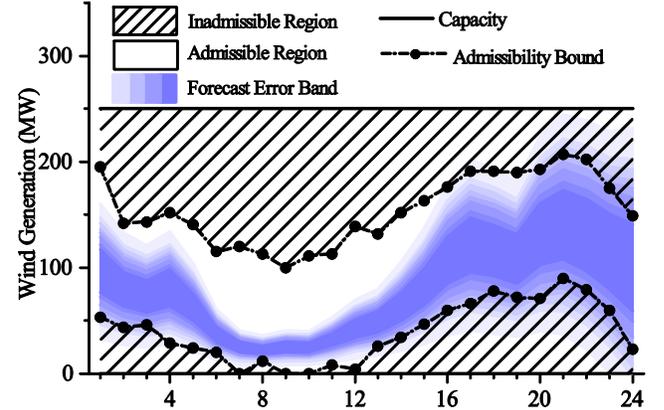

Fig. 1. Schematic diagram of admissible region of wind generation uncertainty.

### B. RRUC Formulation

RRUC aims to minimize the UC and ED cost under the forecasted value of wind generation while guaranteeing the feasibility of dispatch strategy under the adjustable uncertainty set. The RRUC is modeled as problem (2).

In (2), (3a) is to minimize the operational cost, in which the first term represents the UC cost and last two terms represent the ED cost under forecasted wind generation, $\hat{w}_{mt}$. $C_g(\cdot)$ is quadratic and can be further linearized by using PWL method. (3b) and (3c) describe the minimum on/off period limits of generators. (3d) is the start-up constraint of generators. (3e) is the generation capacity of generators. (3f) and (3g) are the ramping rate limits of generators. (3h) depicts the power balance requirement under $\hat{w}_{mt}$. (3i) is the network power flow limit. (3j) is the operational risk limit. (3k) and (3l) depict the boundaries of $w_{mt}^l$ and $w_{mt}^u$, respectively. (1c) and (1d) depict the piecewise linear relationship between $Q_{mt}^p$, $Q_{mt}^n$ and $w_{mt}^u, w_{mt}^l$. $\Omega$ is the feasibility set of $u_{gt}, w_{mt}^u, w_{mt}^l$, which can be defined as problem (3).

$$\min_{z,u,\hat{p},w^u,w^l,Q^p,Q^n} \sum_{t=1}^{T}\sum_{g=1}^{G}\left(S_g z_{gt} + c_g u_{gt} + C_g(\hat{p}_{gt})\right) \quad (3a)$$

$$s.t. \; -u_{g(t-1)} + u_{gt} - u_{gk} \leq 0, \; \forall g,t, k=t,\cdots,t+T_g^{on}-1. \quad (3b)$$

$$u_{g(t-1)} - u_{gt} + u_{gk} \leq 1, \; \forall g,t, k=t,\cdots,t+T_g^{off}-1. \quad (3c)$$

$$-u_{g(t-1)} + u_{gt} - z_{gt} \leq 0, \; \forall g, \forall t \quad (3d)$$

$$u_{gt} P_{\min}^g \leq \hat{p}_{gt} \leq u_{gt} P_{\max}^g \; \forall g, \forall t \quad (3e)$$

$$\hat{p}_{gt} - \hat{p}_{g(t+1)} \leq u_{g(t+1)} R_-^g + (1-u_{g(t+1)}) P_{\max}^g \; \forall g, \forall t \quad (3f)$$

$$\hat{p}_{g(t+1)} - \hat{p}_{gt} \leq u_{gt} R_+^g + (1-u_{gt}) P_{\max}^g \; \forall g, \forall t \quad (3g)$$

$$\sum_{g=1}^{G} \hat{p}_{gt} + \sum_{m=1}^{M} \hat{w}_{mt} = \sum_{j=1}^{J} D_{jt} \; \forall t \quad (3h)$$

$$-F_l \leq \sum_{g=1}^{G} \pi_{gt} \hat{p}_{gt} + \sum_{m=1}^{M} \pi_{mt} \hat{w}_{mt} - \sum_{j=1}^{J} \pi_{jt} D_{jt} \leq F_l \; \forall t \quad (3i)$$

$$Risk = \min_{Q_{mt}^p, Q_{mt}^n} \sum_{t=1}^{T}\sum_{m=1}^{M}(Q_{mt}^p + Q_{mt}^n) \leq Risk_{dh} \quad (3j)$$

$$0 \leq w_{mt}^l \leq \hat{w}_{mt} \; \forall m, \forall t \quad (3k)$$

$$\hat{w}_{mt} \leq w_{mt}^u \leq w_m^{\max} \; \forall m, \forall t \quad (3l)$$

(1c)-(1d)

$$u_{gt}, w_{mt}^u, w_{mt}^l \in \Omega \quad (3m)$$

$$\Omega := \left\{ u_{gt}, w_{mt}^u, w_{mt}^l \,\middle|\, \max_{v^u, v^l} \min_{p, \Delta w, \Delta D} \sum_{t=1}^{T}\left(\sum_{m=1}^{M}\Delta w_{mt} + \sum_{j=1}^{J}\Delta D_{jt}\right) = 0 \quad (3n) \right.$$

$$s.t. \; u_{gt} P_{\min}^g \leq p_{gt} \leq u_{gt} P_{\max}^g \; \forall g, \forall t \quad (3o)$$

$$p_{gt} - p_{g(t+1)} \leq u_{g(t+1)} R_-^g + (1-u_{g(t+1)}) P_{\max}^g \; \forall g, \forall t \quad (3p)$$

$$p_{g(t+1)} - p_{gt} \leq u_{gt} R_+^g + (1-u_{gt}) P_{\max}^g \; \forall g, \forall t \quad (3q)$$

$$\sum_{g=1}^{G} p_{gt} + \sum_{m=1}^{M}(w_{mt} - \Delta w_{mt}) = \sum_{j=1}^{J}(D_{jt} - \Delta D_{jt}) \quad (3r)$$

$$0 \leq \Delta D_{jt} \leq D_{jt} \; \forall j, \forall t \quad (3s)$$

$$0 \leq \Delta w_{mt} \leq w_{mt} \; \forall m, \forall t \quad (3t)$$

$$-F_l \leq \sum_{g=1}^{G}\pi_{gt} p_{gt} + \sum_{m=1}^{M}\pi_{mt}(w_{mt} - \Delta w_{mt}) - \cdots$$
$$- \sum_{j=1}^{J}\pi_{jt}(D_{jt} - \Delta D_{jt}) \leq F_l \; \forall t \quad (3u)$$

$$w_{mt} = (w_{mt}^u - \hat{w}_{mt}) v_{mt}^u + (w_{mt}^l - \hat{w}_{mt}) v_{mt}^l + \hat{w}_{mt} \quad (3v)$$

$$\sum_{t=1}^{T}(v_{mt}^u + v_{mt}^l) \leq \Gamma^T \; \forall m \quad (3w)$$

$$\sum_{m=1}^{M}(v_{mt}^u + v_{mt}^l) \leq \Gamma^S \; \forall t \quad (3x)$$

$$v_{mt}^u + v_{mt}^l \leq 1 \; \forall m, \forall t \quad (3y)$$

$$\left. v_{mt}^u, v_{mt}^l \in \{0,1\} \right\} \quad (3z)$$

(3)

In (3), (3n) is the sum of LS and WGC. (3o) depicts the capacity of generators. (3p) and (3q) limit the ramping capacity of generators. (3r) depicts the relaxed power balance requirement with recourse actions including LS and WGC. (3s) and (3t) are the boundaries of LS and WGC respectively. (3u) is the network power flow limit considering LS and WGC. (3v)-(3z) use a polyhedral set to describe the wind generation denoted by $W$. Specifically, (3v) depicts the wind generation output; (3w) and (3x) describe the uncertainty budgets over both temporal and spatial domains, respectively. Specially, unlike the description of $W$ in the literature [4], [5], the proposed $W$ in (3) is adjustable variables.

In light of problem (2) and (3), RRUC is a two-stage robust optimization problem. The first-stage decision variables are $z_{gt}, u_{gt}, \hat{p}_{gt}, w_{mt}^u, w_{mt}^l, Q_{mt}^p, Q_{mt}^n$; the recourse action variables are $p_{gt}, \Delta w_{mt}, \Delta D_{jt}$; the uncertainty variables are $v_{mt}^u, v_{mt}^l$. Due to the existence of (3n), no LS or WGC will occur in the recourse stage, which guarantees the operational feasibility of $u_{gt}$ as well as the admissibility of $w_{mt}^l$ and $w_{mt}^u$.

In problem (2), noted that (3j) itself is a minimization problem, it can be directly transformed into a standard constraint as follows, as (3a) is also a minimization problem.

$$\sum_{t=1}^{T}\sum_{m=1}^{M}(Q_{mt}^p + Q_{mt}^n) \leq Risk_{dh} \quad (4a)$$

The resulting $w^u$ and $w^l$ are feasible solutions with respect to (3) and (4a), however, may not be the optimal with respect to the resulting UC strategy, as the operational risk is not minimized in the objective function of problem (2), which means risk-based admissibility assessment [24] or DNE limit [15] assessment has to be conducted to obtain the optimal $w^u$ and $w^l$ with respect to the resulting UC strategy. To avoid additional computation burden as well as to remain the simplicity of the proposed framework, one treatment is to add (1b) with penalty coefficient into (3a). Then we have

$$\min_{z,u,\hat{p},w^u,w^l,Q^p,Q^n} \sum_{t=1}^{T}\left(\sum_{g=1}^{G}\left(S_g z_{gt} + c_g u_{gt} + C_g(\hat{p}_{gt})\right) + K \cdot \sum_{m=1}^{M}(Q_{mt}^p + Q_{mt}^n)\right) \quad (4b)$$

With (1b) being minimized in (4b), the resulting $w^u$ and $w^l$ are also optimal with respect to the UC strategy in terms of operational risk minimization. Then RRUC can be reformulated into a standard two-stage robust optimization problem as follows.

*Objective*: (4b)
*s.t.* (1c)-(1d), (3b)-(3i), (3k)-(3m), (4a) (4)

*Remarks:*

1) The constraints of the problem (2) and problem (4) are the same. However, their objective functions are different. Thus, problem (4) is only an approximation of problem (2). The optimality gap between these problems can be controlled by tuning the penalty coefficient $K$. Detailed discussions will be provided in section IV.

2) The value of $Risk_{dh}$ in (4a) is an important parameter that influences not only the RRUC strategy but also the solvability of (4). In practice, it can be chosen depending on historical operation data, operator's risk preference, electricity contract etc.

3) There are many equivalent and tighter forms for some constraints of problem (2), i.e., the minimum on/off time constraints (3b)-(3c) [25]. For tight and strong formulations of UC problem, interested readers can refer to [26], [27].

## C. Compact Model of RRUC

For ease of analysis, the RRUC can be written in a compact form as below:

$$\min_{\mathbf{x},\hat{\mathbf{y}},\mathbf{w},\mathbf{Q}} \mathbf{a}^T\mathbf{x} + \mathbf{b}^T\hat{\mathbf{y}} + \mathbf{c}^T\mathbf{Q} \quad (5a)$$

$$s.t. \quad \mathbf{Ax} + \mathbf{B}\hat{\mathbf{y}} \leq \mathbf{d} \quad (5b)$$

$$\mathbf{Cw} + \mathbf{DQ} \leq \mathbf{e} \quad (5c)$$

$$\begin{pmatrix}\mathbf{x}\\\mathbf{w}\end{pmatrix} \in \begin{cases} \mathbf{x},\mathbf{w} \mid \max_{\mathbf{v}} \min_{\mathbf{y},\mathbf{s}} \mathbf{f}^T\mathbf{s} = 0 & (5d) \\ s.t. \quad \mathbf{Ex} + \mathbf{Fy} + \mathbf{G}(\mathbf{w}\circ\mathbf{v}) + \mathbf{Hs} + \mathbf{Jv} \leq \mathbf{g} & (5e) \\ \mathbf{Lv} \leq \mathbf{h} & (5f) \end{cases}$$

In (5), **x** represents the binary vector of generators. $\hat{\mathbf{y}}$ and **y** represent the continuous vector of generators. **w** represents the wind generation output boundary vector. **Q** represents the operational risk vector. **s** represents the LS as well as WGC vector. **v** is the binary vector depicting wind generation uncertainty. **a, b, c, d, e, f, g, h, A, B, C, D, E, F, G, H, J, L** are constant coefficient matrices and can be derived from (4), in which, i.e., matrix **A** can be derived from (3b)-(3g). '∘' in (5e) is a Hadamard product. Compared with RUC, more variables and constraints are involved in RRUC. Specifically, **w** becomes a decision variable to be determined in the first stage, considerably increasing the computational complexity.

## III. SOLUTION METHODOLOGY

In this section, we will derive the solution algorithms to solve (5). First of all, each stage of (5) is written separately as
*Main Problem (MP):* (5a)-(5c)
*Feasibility & Admissibility Checking Subproblem (F&ACSP):*

$$\max_{\mathbf{v}} \min_{\mathbf{y},\mathbf{s}} \mathbf{f}^T\mathbf{s} \quad (6a)$$
$$s.t. \quad (5e)\text{-}(5f) \quad (6)$$

The solution methodology for F&ACSP is firstly proposed. Then, a C&CG based algorithm is adopted to solve MP. At last, a computational scale reduction algorithm for F&ACSP and a convergence acceleration algorithm for MP are derived.

### A. Solution Methodology for F&ACSP

Mathematically, F&ACSP is a bi-level mixed integer linear program (MILP) and can be solved by many effective methods based on Karush-Kuhn-Tucker (KKT) conditions [28] or the strong duality theory [29]. In this paper, the inner problem of (6a) is replaced by its dual to come up to a single-level bilinear program that can be solved by either the outer approximation method (OA) [16] or the big-M linearization method [30]. As OA may fail to find the global optimal solution in some circumstances, this paper adopts the big-M linearization method to solve (6a) with constraints (5e)-(5f). The compact formulation of dual problem of (6) is as follows. It should be noted that the big-M method may face scalability problem as auxiliary binary variables and constraints would be introduced.

$$\max_{\mathbf{v},\boldsymbol{\lambda}} R = \boldsymbol{\lambda}^T(\mathbf{g} - \mathbf{Ex}) - \boldsymbol{\lambda}^T\mathbf{Jv} - \boldsymbol{\lambda}^T\mathbf{G}(\mathbf{w}\circ\mathbf{v}) \quad (7a)$$

$$s.t. \quad [\mathbf{F} \vdots \mathbf{H}]^T \boldsymbol{\lambda} \leq [\mathbf{0}^T \vdots \mathbf{f}^T]^T \quad (7b) \quad (7)$$

$$\boldsymbol{\lambda} \leq \mathbf{0} \quad (7c)$$

(5f)

where, $\boldsymbol{\lambda}$ is the dual vector of inner problem of (6a). Noting that there are bilinear terms in (7a), auxiliary variables and constraints are introduced to replace them to convert (7) to be a MILP problem as follows.

$$\max_{\mathbf{v},\boldsymbol{\lambda},\boldsymbol{\gamma}} R = \boldsymbol{\lambda}^T(\mathbf{g}\text{-}\mathbf{Ex}) - \boldsymbol{\gamma}^T\mathbf{q} \quad (8a)$$

$$s.t. \quad (5f), (7b)\text{-}(7c)$$
$$-M_{Big}\mathbf{v} \leq \boldsymbol{\gamma} \leq \mathbf{0} \quad (8b) \quad (8)$$
$$-M_{Big}(\mathbf{1}-\mathbf{v}) \leq \boldsymbol{\lambda} - \boldsymbol{\gamma} \leq \mathbf{0} \quad (8c)$$

where, $\boldsymbol{\gamma}$ is the auxiliary vector, **q** is a constant vector and can be derived from the following formula.

$$\boldsymbol{\lambda}^T\mathbf{Jv} + \boldsymbol{\lambda}^T\mathbf{G}(\mathbf{w}\circ\mathbf{v}) = \sum_i\sum_j q_{ij}\lambda_i v_j = \boldsymbol{\gamma}^T\mathbf{q}, \; \gamma_{ij} = \lambda_i v_j \quad (9)$$

(8b) and (8c) are auxiliary constraints generated during objective function linearization using the big-M method. $M_{big}$ is sufficient large positive real number. Thus, (8) result in a standard single-level MILP, which can be solved by using commercial solvers such as CPLEX. From simulation results, solution efficiency of (8) is directly proportional to the scale of $\boldsymbol{\gamma}$ and (8b)-(8c). Meanwhile, it is not difficult to figure out that the scale of $\boldsymbol{\gamma}$ and (8b)-(8c) is related to the number of non-zero elements in **G**. In other words, if the sparsity of **G** can be improved, the computational scale of (8) will be decreased.

### B. Solution Methodology for RRUC Problem

Note that MP (5a) with constraints (5b)-(5c) and F&ACSP (8a) with constraints (5f), (7b)-(7c), (8b)-(8c) both are MILPs. Next the C&CG algorithm is adopted to solve RRUC problem and named as A1. The details of A1 is as follows.

| A1: C&CG Algorithm |
| --- |
| *Step 1: set l=0 and $\mathbf{O} = \emptyset$.* |
| *Step 2: Solve (5a)-(5c) with the additional constraints as follows.* |
| $\mathbf{Ex} + \mathbf{Fy}^k + \mathbf{G}(\mathbf{w}\circ\mathbf{v}_k^*) + \mathbf{Jv}_k^* \leq \mathbf{g} \; \forall k \leq l$ *(10a)* |
| *Step 3: Solve model (8). If $\|R_{k+1} - R_k\| < \epsilon$, terminate. Otherwise, derive the optimal solution $\mathbf{v}_{k+1}^*$, create variable vector $\mathbf{y}^{k+1}$ and add the following constraints* |
| $\mathbf{Ex} + \mathbf{Fy}^{k+1} + \mathbf{G}(\mathbf{w}\circ\mathbf{v}_{k+1}^*) + \mathbf{Jv}_{k+1}^* \leq \mathbf{g}$ *(10b)* |
| *Update l=l+1, $\mathbf{O} = \mathbf{O}\cup\{l+1\}$ and go to Step 2.* |

In A1, $\epsilon$ represents the convergence gap. In the standard C&CG algorithm [31], a set of constraints (5e) of F&ACSP with the identified worst-case scenario are directly added into MP. However, in A1, the added constraints (10b) are not the same with the original constraints (5e) in F&ASP. Compared with (10b), (5e) can be regarded as loose constraints with slack variables as recourse actions are involved.

### C. Computational Scale Reduction

As mentioned above, the key point to improve the efficiency in solving (8) is to enhance the sparsity of **G**. To do this, we replace (3r) and (3u) with the following constraints.

$$\sum_{g\in\phi(n)}p_{gt} + \sum_{m\in\phi(n)}(w_{mt} - \Delta w_{mt}) - \sum_{o\in\phi(n)}B_{on}(\theta_{nt} - \theta_{ot}) \\ - \sum_{j\in\phi(n)}(D_{jt} - \Delta D_{jt}) = 0 \quad \forall n, \forall t \quad (10c)$$

$$-F_l \le B_{o_1 o_2}(\theta_{o_1 t} - \theta_{o_2 t}) \le F_l \quad o_1, o_2 \in Line_l, \forall l, \forall t \quad (10d)$$

$$-\pi \le \theta_{nt} \le \pi \quad \forall n, \forall t \quad (10e)$$

$$\theta_{reft} = 0 \quad \forall t \quad (10f)$$

Constraint (10c) represents the power balance equation for each node. (10d) is the power flow limit on transmission lines. (10e) describes the upper and lower limits of the nodal phase angles and (10f) represents the reference phase angle. In other words, network power balance constraint (3r) is replaced by nodal power balance constraint (10c). Moreover, transmission limit (3u) based on nodal power injection sensitivity matrix (NPISM) is replaced by (10d) based on phase angle and node admittance matrix. Similarly, the compact formulation of F&ACSP with replaced constraints is as follows.

$$\max_{\mathbf{v}} \min_{\mathbf{z},\mathbf{s}} \mathbf{e}^T \mathbf{s} \quad (11a)$$

$$s.t. \quad \mathbf{Mx} + \mathbf{Nz} + \mathbf{O}(\mathbf{w} \circ \mathbf{v}) + \mathbf{Ps} + \mathbf{Uv} \le \mathbf{p} \quad (11b) \quad (12)$$

(5f)

In (13), $\mathbf{z}$ is the continuous vector including output of generators and phase angle of each node. $\mathbf{M}, \mathbf{N}, \mathbf{O}, \mathbf{P}, \mathbf{U}, \mathbf{e}, \mathbf{p}$ are constant coefficient matrices, which can be derived from (3n)-(3q), (3s)-(3z) and (10c)-(10f), respectively. Similarly, (12) can be rewritten as a single level linear problem as follows.

$$\max_{\mathbf{v},\boldsymbol{\eta},\boldsymbol{\mu}} R = \boldsymbol{\eta}^T (\mathbf{p} - \mathbf{Mx}) - \boldsymbol{\mu}^T \mathbf{r} \quad (12a)$$

$$s.t. \quad [\mathbf{N} \vdots \mathbf{P}]^T \boldsymbol{\eta} \le [\mathbf{0}^T \vdots \mathbf{e}^T]^T \quad (12b)$$

$$\boldsymbol{\eta} \le \mathbf{0} \quad (12c)$$

$$-M_{Big} \mathbf{v} \le \boldsymbol{\mu} \le \mathbf{0} \quad (12d) \quad (13)$$

$$-M_{Big}(\mathbf{1} - \mathbf{v}) \le \boldsymbol{\eta} - \boldsymbol{\mu} \le \mathbf{0} \quad (12e)$$

(5f)

In (13), $\boldsymbol{\eta}$ is the dual vector and $\boldsymbol{\mu}$ is the auxiliary vector, $\mathbf{r}$ is a constant vector and can be derived from the following formula.

$$\boldsymbol{\eta}^T \mathbf{Uv} + \boldsymbol{\eta}^T \mathbf{O}(\mathbf{w} \circ \mathbf{v}) = \sum_i \sum_j r_{ij} \eta_i v_j = \boldsymbol{\mu}^T \mathbf{r}, \quad \mu_{ij} = \eta_i v_j \quad (13)$$

Generally, the number of non-zero elements in $\mathbf{O}$ is much smaller than that in $\mathbf{G}$. Comparison of computational complexity between (8) and (13) is listed in Table I. From Table I, although the number of continuous variables in (13) is larger than that in (8) by $(3N+1)T$ and the number of regular constraints with respect to $\boldsymbol{\eta}, \mathbf{v}$ in (13) is larger than regular constraints with respect to $\boldsymbol{\lambda}, \mathbf{v}$ in (8) by $NT$, the numbers of the rest variables and constraints in (13) are much smaller than that in (8), especially in the big-M constraints. Here the second algorithm is derived and named as A2. The only difference from A2 is that (13) instead of (8) is solved in step 3 of A2. For simplicity, details of A2 are omitted.

TABLE I COMPUTATIONAL SCALE COMPARISON

|   | Model (8) | Model (13) |
|---|---|---|
| Binary Variables | $\mathbf{v}$: $2MT$ | $\mathbf{v}$: $2MT$ |
| Continuous Variables | $\boldsymbol{\lambda}$: $(3G+2L+2J+2M)T$ | $\boldsymbol{\eta}$: $(3G+2L+2J+2M+\cdots$ $\cdots+3N+1)T$ |
| Auxiliary Variable | $\boldsymbol{\gamma}$: $4(L+1)MT$ | $\boldsymbol{\mu}$: $4MT$ |
| Regular Constraints | $\boldsymbol{\lambda}, \mathbf{v}$: $(G+M+L)T$ | $\boldsymbol{\eta}, \mathbf{v}$: $(G+M+L+N)T$ |
| Regular Constraints | $\boldsymbol{\lambda}, \boldsymbol{\gamma}$: $8(L+1)MT$ | $\boldsymbol{\lambda}, \boldsymbol{\mu}$: $8MT$ |
| Big-M Constraints | $\boldsymbol{\lambda}, \mathbf{v}, \boldsymbol{\gamma}$: $8(L+1)MT$ | $\boldsymbol{\eta}, \mathbf{v}, \boldsymbol{\mu}$: $8MT$ |

Compared with the standard C&CG algorithm, on one hand, A2 increases the computational scale of inner problem of (6), which leads to the decrement of computational burden of F&ACSP in return; on the other hand, the scale of variables and constraints generated in each iteration of MP are remarkably decreased compared with passing the variables and constraints of F&ACSP to MP directly.

It should be pointed out that the computational scale reduction approach discussed above is also applicable to other two-stage robust optimization problems in which power balance constraint is involved in the second stage such as RUC. Further, the effectiveness of this approach will be more significant with the increase of wind farms.

*D. Convergence Acceleration*

In some cases, the convergence efficiency of A2 is not very satisfying. In this regard, some active constraints are generated and added into MP in each iteration to speed up the convergence. In light of [16], a feasibility cut is generated in each iteration and added into MP. The construction of feasibility cut in iteration $k+1$ is

$$-\boldsymbol{\eta}_k^T \mathbf{M}(\mathbf{x} - \mathbf{x}_k^*) - \boldsymbol{\eta}_k^T \mathbf{O}(\mathbf{w} \circ \mathbf{v}_k^* - \mathbf{w}_k^* \circ \mathbf{v}_k^*) \le -R_k \quad (14)$$

In (14), $\mathbf{x}_k^*, \mathbf{w}_k^*$ are the optimal solutions of MP in iteration $k$. $\eta_k$ is the optimal solution of F&ACSP in iteration $k$. $R_k$ is the objective value of MP in iteration $k$. Actually, (14) serves as the sub-gradient cut and gradient cut for $\mathbf{x}$ and $\mathbf{w}$, respectively. The third algorithm is developed as follows and named as A3.

| A3: C&CG-based Algorithm with feasibility cut |
|---|
| *Step 1:* set $l=0$ and $\mathbf{O} = \emptyset$. |
| *Step 2:* Solve (5a)-(5c) with the additional constraints as follows. |
| $\quad \mathbf{Ex} + \mathbf{Fy}^k + \mathbf{G}(\mathbf{w} \circ \mathbf{v}_k^*) + \mathbf{Jv}_k^* \le \mathbf{g} \quad \forall k \le l \quad (15a)$ |
| $\quad -\boldsymbol{\eta}_k^T \mathbf{M}(\mathbf{x} - \mathbf{x}_k^*) - \boldsymbol{\eta}_k^T \mathbf{O}(\mathbf{w} \circ \mathbf{v}_k^* - \mathbf{w}_k^* \circ \mathbf{v}_k^*) \le -R_k \quad \forall k \le l \quad (15b)$ |
| *Step 3:* Solve model (14). If If $|R_{k+1} - R_k| < \epsilon$, terminate. Otherwise, derive the optimal solution $\mathbf{v}_{k+1}^*, \eta_{k+1}$, create variable vector $\mathbf{y}^{k+1}$ and augment the following constraints |
| $\quad \mathbf{Ex} + \mathbf{Fy}^{k+1} + \mathbf{G}(\mathbf{w} \circ \mathbf{v}_{k+1}^*) + \mathbf{Jv}_{k+1}^* \le \mathbf{g} \quad (15c)$ |
| $\quad -\boldsymbol{\eta}_{k+1}^T \mathbf{M}(\mathbf{x} - \mathbf{x}_{k+1}^*) - \boldsymbol{\eta}_{k+1}^T \mathbf{O}(\mathbf{w} \circ \mathbf{v}_{k+1}^* - \mathbf{w}_{k+1}^* \circ \mathbf{v}_{k+1}^*) \le -R_{k+1} \quad (15d)$ |
| *Update* $l=l+1$, $\mathbf{O} = \mathbf{O} \cup \{l+1\}$ and go to Step 2. |

Compared with A2, the feasibility cut (14) and the value of dual variable $\boldsymbol{\eta}$ of F&ACSP are passed to MP in each iteration in A3. The computational efficiency of A1-A3 will be compared in Section IV. D.

## IV. CASE STUDIES

In this section, numerical experiments on the modified IEEE 118-bus system are carried out to show the effectiveness of the proposed model and algorithms. The experiments are performed on a PC with Intel(R) Core(TM) 2 Duo 2.2 GHz CPU and 4 GB memory. All algorithms are implemented on MATLAB and programmed using YALMIP. The MILP solver is CPLEX 12.6. The optimality gap is set as 0.1%.

*A. The Modified IEEE 118-bus System*

The tested system has 54 generators and 186 transmission lines. Three wind farms are connected to the system at bus 17, 66 and 94, respectively. The installed capacities are 500 MW identically. The generators' parameters and the load curve can be found in [32]. All day-ahead forecast of wind generation are scaled down from the day-ahead curve of California ISO as shown in Fig. 2. Prices for LS and WGC are listed in Table II. We choose the confidence level $\beta_t = 95\%$ and $\beta_s = 95\%$, yielding $\Gamma^T \approx 8$ and $\Gamma^S \approx 2$ [16]. In this case, the root mean square errors of $\delta_{mt}$ are subject to (15) with $\sigma_1 = 20\%$, $\sigma_2 = 15\%$, $\sigma_3 = 10\%$ and their mean value are zero. In (15), $\sigma_m$ is a constant parameter.

$$\sigma_{mt} = \sigma_m \cdot \hat{w}_{mt} \cdot (1 + e^{-(T-t)}) \quad \forall m, \forall t \tag{15}$$

Wind generation forecast error bands are simply derived by Gaussian distribution and as shown in Fig. 3. There are other advanced methods to determine wind generation forecast error bands in the literature. However, it has no influence on the computation solvability and is beyond the scope of this paper. We choose $\alpha_1^n = 0.5\%, \alpha_2^n = 2.5\%, \alpha_3^n = 49.5\%$, which means an eight-piecewise linear distribution approximation is adopted to depict the PDF of $\delta_{mt}$. Details of the PDF approximation can be found in [24]. We further set $Z=4$ in (1c) and (1d) based on the setting presented in [24].

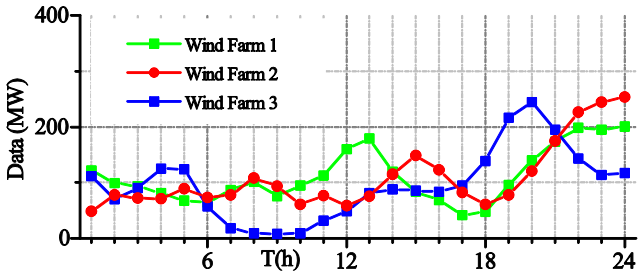
Fig. 2. Forecasted value of wind farm 1-3.

TABLE II. COST COEFFICIENT OF LS AND WGC IN DIFFERENT TIME PERIODS

| Period | T1:1-6 | T2:7-12 | T3:13-18 | T4:19-24 |
|---|---|---|---|---|
| LS ($/ MWh) | 100 | 200 | 150 | 200 |
| WGC ($/ MWh) | 20 | 40 | 30 | 40 |

### B. Comparison with Other UC Model

In this subsection, RRUC are compared with other UC models in terms of operational cost, operational risk and operational loss, respectively. Specifically, the deterministic unit commitment (DUC) model is from [16], in which the spinning reserve rate is 10%; the SUC model as well as the scenario generation and reduction method are from [2], in which 200 scenarios are originally generated and 20 scenarios are left after the reduction; the RUC model is from [16], in which confidence level $\alpha_t$ is chosen as 95%. The operational risk of RUC is evaluated based on the method presented in [24]. Then the evaluated operational risk is regarded as the benchmark and is selected as $Risk_{dh}$ for RRUC. According to the value of $Risk_{dh}$, $K$ is selected as 0.1. The operational cost, operational risk as well as computational time under those four UC models are listed in Table III. From Table III, both operational cost and operational risk of RRUC are lower than RUC, showing a better capability of optimizing operational flexibility as well as mitigating operational risk. Meanwhile, the sum of operational cost and operational risk of RRUC are the lowest among those UC models. It should be noted that, the computational time of RRUC is the highest among the four UC models, if the scenario reduction time of SUC is not considered. Detailed computation efficiency analysis will be demonstrated in Section IV.D.

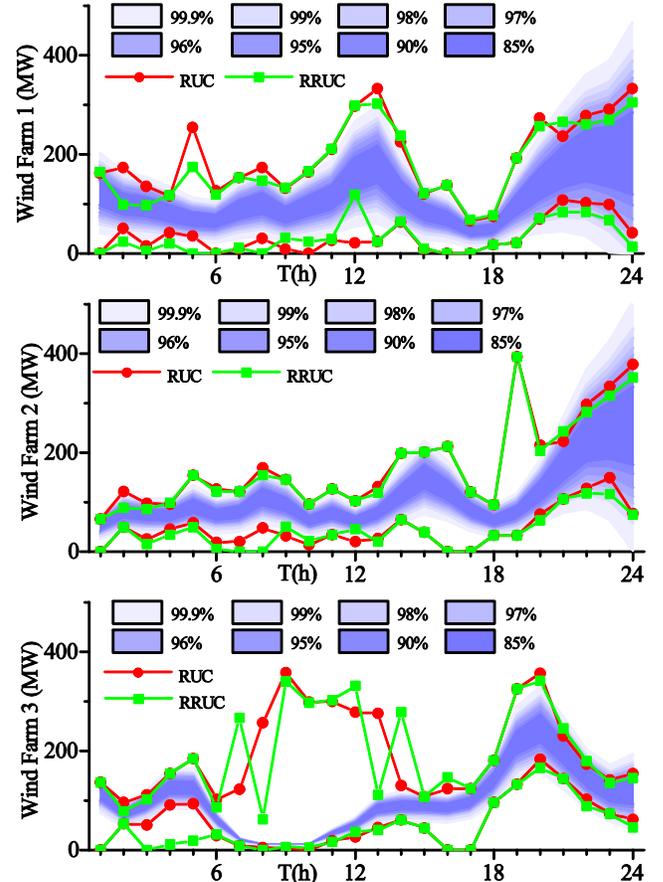
Fig. 3. Admissible wind generation boundary of different UC strategies.

TABLE III COST AND RISK UNDER DIFFERENT UC MODELS

| | Total Cost ($) | UC Cost ($) | ED Cost ($) | Risk ($) | Time (s) |
|---|---|---|---|---|---|
| DUC | 1.287×10^6 | 1.90×10^4 | 1.262×10^6 | 2.67×10^4 | 125 |
| SUC | 1.304×10^6 | 2.79×10^4 | 1.276×10^6 | 9.86×10^3 | 1879 |
| RUC | 1.312×10^6 | 3.29×10^4 | 1.283×10^6 | 7.23×10^3 | 3727 |
| RRUC | 1.307×10^6 | 2.87×10^4 | 1.278×10^6 | 6.64×10^3 | 5399 |

TABLE IV OPERATIONAL LOSS OF DIFFERENT UCS UNDER RARE EVENTS

| | Average Operational Loss ($) | | |
|---|---|---|---|
| | Total | WGC | LS |
| DUC | 1.017×10^6 | 2.172×10^5 | 8.010×10^5 |
| SUC | 5.094×10^5 | 1.365×10^5 | 3.720×10^5 |
| RUC | 4.050×10^5 | 1.209×10^5 | 2.841×10^5 |
| RRUC | 3.357×10^5 | 1.317×10^5 | 2.043×10^5 |

We define a wind generation scenario being partly or fully out of a prescribed wind generation uncertainty set as a rare event. To test the performance of those UCs under rare events, 10,000 wind generation scenarios are generated beyond the wind generation uncertainty set ($\alpha_t = 95\%$) of RUC. The results are demonstrated in Table IV. From Table IV, RRUC gives the lowest total average operational loss, which confirms the operational risk index in Table III. Also, the wind generation admissibility boundaries under RUC and RRUC are given in Fig. 3. From Fig. 3, both the upper and the lower admissible boundaries of RRUC are lower than RUC in most periods. This explains the differences of operational loss resulted from WGC and LS between RUC and RRUC.

### C. Uncertainty Set and Admissibility Region

In RUC, the uncertainty set is prescribed, which is equal to the 95% forecast error band and denoted as $W^{RUC}$. Therefore, $W^{RUC}$ is symmetric with respect to $\hat{w}_{mt}$ and the width of $W^{RUC}$ in each period is proportional to $\sigma_{mt}$. The admissibility region of RUC, however, is always larger than $W^{RUC}$, denoted as $R^{RUC}$ (red lines in Fig. 3). In RRUC, the uncertainty set, denoted as $W^{RRUC}$ (green lines in Fig. 3), is variable and its width varies in both spatial and temporal domains, which reflects the optimal allocation of operational flexibility as well as operational risk mitigation capability, resulting in operational cost and risk decrement than RUC. Besides, the admissibility region of RRUC, denoted as $R^{RRUC}$ (green lines in Fig. 3), is identical with $W^{RRUC}$.

### D. Computational Efficiency

In this subsection, the computational efficiency of A1-A3 under different uncertainty budget $\varGamma^T$ are discussed. The simulation results are listed in Table V. It is observed that A2 enhances the computational efficiency by 83.3% averagely compared with A1 due to the computational scale reduction. A3 further improves the computational efficiency by 77.5% in average compared with A2, by 225% in average compared with A1 by decreasing iteration number and reducing computational scale simultaneously. These results manifest the effectiveness of A2 and A3.

TABLE V
COMPUTATIONAL EFFICIENCY UNDER DIFFERENT CASES AND ALGORITHMS

|    |              | Total (s) | MP (s) | F&ACSP (s) | Iteration |
|----|--------------|-----------|--------|------------|-----------|
| A1 | $\varGamma^T$=8  | 15892 | 8971 | 6921 | 15 |
|    | $\varGamma^T$=16 | 7239  | 3406 | 3833 | 9  |
|    | $\varGamma^T$=24 | 3753  | 1082 | 2671 | 5  |
| A2 | $\varGamma^T$=8  | 9775  | 8614 | 1161 | 15 |
|    | $\varGamma^T$=16 | 3647  | 3013 | 634  | 9  |
|    | $\varGamma^T$=24 | 1255  | 992  | 263  | 5  |
| A3 | $\varGamma^T$=8  | 5399  | 4587 | 812  | 12 |
|    | $\varGamma^T$=16 | 2183  | 1811 | 372  | 7  |
|    | $\varGamma^T$=24 | 691   | 590  | 101  | 4  |

TABLE VI SIMULATION RESULTS UNDER DIFFERENT VALUE OF K

|        | Total Cost ($) | UC Cost ($) | ED Cost ($) | Risk($) |
|--------|----------------|-------------|-------------|---------|
| K=0.1  | 1.3067×10$^6$ | 2.874×10$^4$ | 1.278×10$^6$ | 6.64×10$^3$ |
| K=1    | 1.3067×10$^6$ | 2.874×10$^4$ | 1.278×10$^6$ | 6.64×10$^3$ |
| K=10   | 1.3086×10$^6$ | 2.970×10$^4$ | 1.279×10$^6$ | 6.38×10$^3$ |

### E. Impact of penalty coefficient K

As stated before, the intention adding $K$ into (4b) is to control gap between optimal value of problem (2) and (4). One treatment is to use adaptive $K$ to make the order of magnitude of penalty term lower than that of precision tolerance of MP. In this case, the optimality gap of MP is 0.1% and the order of magnitude of optimal value of MP is $10^6$. Then the order of magnitude of precision of MP is $10^3$. Meanwhile, the order of magnitude of $Risk_{dh}$ is $10^3$ in this case, therefore $K$ can be selected as 0.1 to make the order of magnitude of penalty to be no more than $10^2$, which decreases the gap between these two problems at the optimal solution. Simulation results under different value of $K$ are listed in Table VI. From Table VI, the optimal value of (4) remains unchanged while $K$ varies from 0.1 to 1, showing the effectiveness of proposed method to choose the value of $K$.

### F. Impact of Risk Level

Operational cost and risk of RRUC under different risk levels are shown in Fig. 4. Along with the decrease of the risk level, the operational cost increases gradually. However, when the risk level decreases to a certain critical value, 270$ in this case, RRUC will have no solution anymore if risk level continues to decrease. It means the minimum feasible risk level (MFRL) in this case is 270$. Due to the noncontinuity of UC variables, there exists certain gap between operational risk and risk level. One observation is that the relationship between operational risk and risk level is not strictly linear, as shown in Fig. 4. Similarly, upper bound (UB) and lower bound (LB) of operational risk can also be obtained while changing risk level.

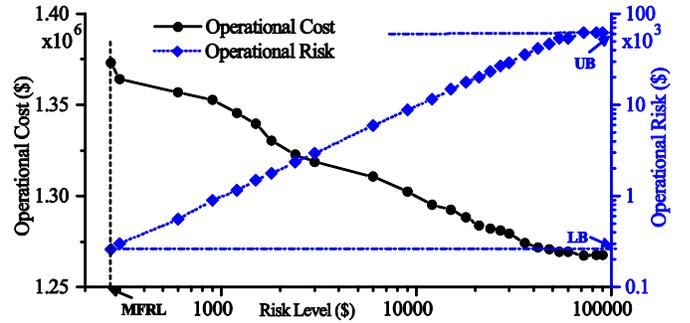

Fig. 4. Operational cost of RRUC under different operational risk levels.

### G. Impact of Numbers of Wind Farms

In this subsection, the impact of numbers of wind farms is analyzed. Here we divide each of the aforementioned wind farms into 2, 3 and 4 small wind farms equally. Computational time under different numbers of wind farms are listed in Table VII, in which $\varGamma^T$=24. From Table VII, as the number of wind farms increases, the computational time increases rapidly, especially the solution time of F&ACSP. It is obvious that if the number of wind farms keep increasing, the big-M based method may not be suitable to solve F&ACSP in terms of computation burden. In this regard, the aforementioned OA approach may be an alternative.

TABLE VII COMPUTATIONAL PERFORMANCE UNDER DIFFERENT WIND FARMS

| Wind Farm | Budget | Total (s) | MP (s) | F&ACSP (s) | Iteration |
|-----------|--------|-----------|--------|------------|-----------|
| 6  | $\varGamma^S$=4 | 748  | 601 | 147 | 4 |
| 9  | $\varGamma^S$=6 | 876  | 627 | 249 | 4 |
| 12 | $\varGamma^S$=8 | 1736 | 891 | 845 | 5 |

## V. CONCLUSION

In this paper, a RRUC model is proposed for the purpose of determining the optimal day-ahead UC strategy under a certain level of operational risk. In the proposed formulation, RRUC is formulated as a two-stage robust optimization problem in which a piecewise linear relationship is constructed between the boundaries of wind generation uncertainty set and operational risk. Compared with RUC, the boundaries of wind generation uncertainty set are adjustable variables to be optimized in RRUC, resulting in an optimal allocation of operational flexibility as well as operational risk mitigating capability.

The proposed methodology answers two important issues missing in the literature of robust UC decision making. on one hand, introducing the concept of operational risk and utilizing it as a constraint solve the problems that how to measure the potential loss if the realization of uncertainty is beyond the prescribed uncertainty set, and how to control it within an acceptable level; on the other hand, the solution inherently provides the prescribed uncertainty set with an optimal admissible boundary in the sense of minimum operational risk.

Three iterative algorithms are proposed based on the C&CG algorithm to solve RRUC effectively. Simulations are carried out on the modified IEEE 118-bus system to illustrate the effectiveness of the proposed model and algorithms. It also reveals the influence of risk level on RRUC.

Due to the limitation of two-stage modeling framework of RRUC, it may not be able to be directly applied to power system decision making, as the operation stage of practical power markets are multiple. However, RRUC can still be applied to part of power market decision making, such as RAA or RAC, considering its modeling flexibility.